\documentclass[journal,twoside,web]{ieeecolor}
\usepackage{lcsys}
\usepackage{graphicx}
\usepackage{graphics}
\usepackage{textcomp}
\usepackage{epsfig}
\usepackage{times}
\usepackage{amsmath}
\usepackage{amssymb}
\usepackage{amsfonts}
\usepackage{float}
\usepackage{siunitx}
\usepackage{scalerel}
\usepackage{cite}
\usepackage{textcomp}
\usepackage{algorithm}
\usepackage{algorithmic}
\usepackage{stackengine}
\usepackage{mathtools}
\def\BibTeX{{\rm B\kern-.05em{\sc i\kern-.025em b}\kern-.08em
    T\kern-.1667em\lower.7ex\hbox{E}\kern-.125emX}}
\def\delequal{\mathrel{\ensurestackMath{\stackon[1pt]{=}{\scriptstyle\Delta}}}}
\newtheorem{theorem}{Theorem}

\newtheorem{lemma}{Lemma}

\newtheorem{remark}{Remark}
 \newtheorem{assumption}{Assumption}
 \newtheorem{Corollary}{Corollary}

\title{\LARGE \bf Distributed Estimation}
\newcommand{\ncom}{\newcommand}

\newcommand{\beqn}{\begin{eqnarray*}}
\newcommand{\eeqn}{\end{eqnarray*}}
\newcommand{\beq}{\begin{eqnarray}}
\newcommand{\eeq}{\end{eqnarray}}

\newcommand{\norm}[1]{\left\lVert #1 \right\rVert}
\newcommand{\inprod}[2]{\left\langle #1, #2 \right\rangle}
\ncom\R{\mathbb{R}}
\DeclareMathOperator*{\argmin}{arg\,min}

\author{Anik Kumar Paul, Arun D Mahindrakar and Rachel K Kalaimani
\thanks{Anik  is a Graduate student in the Department of Electrical Engineering, IIT Madras, Chennai-600036, India
        {email: anikpaul42@gmail.com}}
        \thanks{Arun  and Rachel are  with the Department of Electrical Engineering, Indian Institute of Technology Madras, Chennai-600036, India (email: arun\_dm@iitm.ac.in, rachel@ee.iitm.ac.in) }
    }

\begin{document}
\title{Robust Analysis of Almost Sure Convergence of Zeroth-Order Mirror Descent Algorithm}
\maketitle
\thispagestyle{plain}
\pagestyle{plain}

\begin{abstract}
This letter presents an almost sure convergence of the zeroth-order mirror descent algorithm. The algorithm admits non-smooth convex functions and a biased oracle which only provides noisy function value at any desired point. We approximate the subgradient of the objective function using Nesterov's Gaussian Approximation (NGA) with certain alternations suggested by some practical applications. We prove an almost sure convergence of the iterates' function value to the neighbourhood of optimal function value, which can not be made arbitrarily small, a manifestation of a biased oracle. This letter ends with a concentration inequality, which is  a finite time analysis that predicts the likelihood that the function value of the iterates is in the neighbourhood of the optimal value at any finite iteration.
\end{abstract}


\section{Introduction}
One of the earliest subfields of optimization is derivative-free optimization \cite{10.1093/comjnl/3.3.175,10.1214/aoms/1177699070,brent2002algorithms}. It refers to an optimization problem with an oracle that only provides noisy function value at a desired point. Following numerous attempts by researchers to accurately approximate a function's subgradient from its value (for example see, \cite{doi:https://doi.org/10.1002/0471722138.ch5,NestSpok17}), it has now gained popularity in the optimization community due to its use in a variety of different domains. For a full introduction to derivative-free optimization and its various applications in diverse domains, see \cite{9186148} (and the references therein).

In this letter, we focus on the zeroth-order mirror descent algorithm \cite{9416872}, where the approximated subgradient established in \cite{NestSpok17} replaces the subgradient of the convex objective function in standard mirror descent algorithm \cite{10.1561/2200000050}.  Originally, the mirror descent algorithm generalizes the standard gradient descent algorithm in a more general non-Euclidean   space \cite{anik}.   In recent years, the mirror descent algorithm has grasped significant attention in the large-scale optimization problems, data-driven control and learning, power system, robotics and game theoretic problems \cite{8409957}. For the stochastic mirror descent algorithm, we refer the reader to \cite{doi:10.1137/120894464}.  However, precise information regarding the convex objective function's subgradient or stochastic subgradient is accessible in these articles. In this letter, we assume that we can only access the noisy evaluation of the convex objective function at a desired point via a ``biased zeroth-order" oracle.  The oracle setting is driven by a large number of practical applications in which only the noisy function values are provided at a point and obtaining a subgradient or stochastic subgradient may not be feasible at that point.  As a result, we must approximate the function's subgradient from the noisy measurement of the function value. This gives rise to  the notion of zeroth-order optimization  \cite{7055287}. Every step  in the zeroth-order algorithm is similar to its first-order counterpart (such as gradient descent or mirror descent), except that the function's subgradient must be approximated at every point. There has recently been a surge of interest generated in different variants of  zeroth-order optimization, for both convex and non-convex functions \cite{doi:10.1137/120880811, https://doi.org/10.48550/arxiv.1809.06474,6870494,8619028,8703066, doi:10.1137/18M119046X}, where the subgradient is approximated by NGA \cite{NestSpok17}.

We extend the analysis of zeroth-order optimization in this letter, focusing on the zeroth-order mirror descent (ZOMD) algorithm. The problem framework and analysis in this work differ significantly from the recent literature. The main objective of this study is to show the almost sure convergence of the function value of  iterates of the ZOMD algorithm to a neighbourhood of optimal value, as compared to the bulk of the literature, which focuses on showing that the expected error in function value converges to the neighbourhood of zero. An almost sure convergence guarantee to a neighbourhood of optimal value is more significant than the convergence in expectation since it  describes what happens to the individual trajectory in each iteration.
To the best of our knowledge, no prior work on  almost sure convergence  for zeroth-order optimization has been published.   The problem framework in this study differs from most other works in that it includes a biased oracle that delivers only biased measurement of function value (the expectation of noise in the function measurement is non-zero) at any specified point. The motivation to consider ``biased oracle" can be found in application of reinforcement learning and financial risk measurement (see \cite{prasanth9736646} and references therein for more details).
Furthermore, unlike other publications, we consider that the oracle returns distinct noise values for two different points.
Lastly, in addition to showing almost sure convergence, we estimate the likelihood that the function value of the iterates will be in the neighbourhood of optimal value in any finite iteration. This analysis aids in determining the relationship between the convergence of the ZOMD algorithm and the various parameters of the approximated subgradient. The following list summarises the key contribution of this study.
\begin{enumerate}
    \item  We analyse the ZOMD algorithm under the assumption that a biased oracle returns noisy function value at a predetermined point where the expected error is nonzero. For the biased oracles, we re-evaluate the parameters of the approximated subgradient of the objective function at a specific location, which is calculated using NGA.
    \item We prove that, under certain assumptions, the function values of the iterates of ZOMD algorithm almost surely converges to the neighbourhood of optimal function value. This neighbourhood is determined by several parameters, which are explored in this study.
    \item Finally, we show that for any confidence level and a given neighbourhood around the optimal function value, the function value of the iterate sequence should be in that neighbourhood after some finite iteration with that confidence. We also present an expression for that finite iteration that is influenced by the neighbourhood, confidence level, and other properties stated in the letter.
\end{enumerate}
\section{Notation and Mathematical Preliminaries}
Let $\R$ and $\mathbb{R}^n$ represent the set of real numbers, set of $n$ dimensional real vectors.   Let~$\norm{.}$ denote any \mbox{norm} on~$\R^n$.  Given a norm $\norm{.}$ on $\mathbb{R}^n$, the dual norm of $x\in \R$ is $\norm{x}_\ast := \sup\{ \inprod{x}{y}: {\norm{y}\leq 1}, y\in \R^n \}$, where $\inprod{x}{y}$ denotes the standard inner-product on $\R^n$. $I_n$ is $n \times n$ identity matrix. A random vector  $X\sim \mathcal{N}(0_n, I_n)$ denotes a $n$- dimensional normal random vector with zero-mean and unit standard-deviation. For two random variables $X$ and $Y$,  $\sigma(X,Y)$ is the smallest sigma-algebra generated by  random variables $X$ and $Y$. Because of equivalence of norm $\norm{.}_2 \leq \kappa_1 \norm{.}_\ast$ and $\norm{.}_2 \leq \kappa_2 \norm{.}$ and $\kappa = \kappa_1  \kappa_2$.

Let  $f: \mathbb{R}^n \to \mathbb{R}$ be a convex function.   For $\delta \geq 0$, the vector $g_{\delta}\in \R^n$  is called a $\delta$-subgradient of $f$ at $x$ if and only if $f(y)\geq f(x)+\inprod{g_{\delta}}{y-x} - \delta \;\ \; \forall \; y\in \R^n $ \cite{NestSpok17}. The set of all $\delta$-subgradients at a point $x$ is called the $\delta$-subdifferential of  $f$, denoted by $\partial_\delta f(x)$. If $\delta = 0$, we simply write the notation $\partial f(x)$. If $f$ is differentiable at $x$, then $\partial f(x) = \{ \nabla f(x) \}$, gradient of $f$ at $x$. We say $f \in \mathcal{C}^{0,0}$ if $\exists \; L_0 > 0$ such that $\norm{f(x)-f(y)} \leq L_0 \norm{x-y}$ and $f \in \mathcal{C}^{1,1}$ if $f$ is continuously differentiable and  $\exists \; L_1 >0
$ such that $\norm{\nabla f(x) - \nabla f(y)} \leq L_1 \norm{x-y}$ $\forall \; \; x , y \in \mathbb{R}^n$.

If $f$ has directional derivative in all directions, then we can form the Gaussian approximation as follows: $f_\mu (x) = \frac{1}{(2\pi)^{\frac{n}{2}}}  \int\limits_{\mathbb{R}^n} f(x+\mu u) e^{-\frac{1}{2}\norm{u}^2} du$, where $\mu >0$ is any constant. The function $f_
\mu$ is differentiable at each $x \in \mathbb{R}^n$ and  $\nabla f_{\mu} (x) = \frac{1}{(2 \pi)^{\frac{n}{2}}\mu} \int\limits_{\mathbb{R}^n} u f(x+ \mu u) e^{-\frac{1}{2}\norm{u}^2} du$. It can also be seen that $\nabla f_\mu(x)$ $\in$ $\partial_{\delta} f(x)$, where, $\delta = \mu L_0 \sqrt{n}$ if $f \in \mathcal{C}^{0,0}$ and $\delta = \frac{\mu^2}{2}L_1 \sqrt{n}$ if $f \in \mathcal{C}^{1,1}$.

 Let $(\Omega, \mathcal{F}, \mathbb{P})$ denote a  probability space. An event $A \in \mathcal{F}$ is occurred almost surely (a.s.) if $\mathbb{P}(A) =1$. If $X\sim \mathcal{N}(0_n, I_n)$, it can be shown that $\mathbb{E}[\norm{X}_2^p] \leq n^{\frac{p}{2}}$ if $p \in [0,2]$ and $\mathbb{E}[\norm{X}_2^p] \leq (p+n)^{\frac{p}{2}}$ if $ p > 2$.
We will use the following two Lemmas in our analysis.
\begin{lemma}[  \cite{ROBBINS1971233}]
Let $\{X_t\}_{t \geq 1}$ be a martingale with respect to a filtration $\{\mathcal{F}_t\}_{t \geq 1}$  such that $\mathbb{E}[\norm{X_t}] < \infty$
and $\{\beta(t)\}$ is a non-decreasing sequence of positive numbers such that $\lim\limits_{t \to \infty} \beta(t) = \infty$ \\and $\sum\limits_{t \geq 1} \frac{\mathbb{E}[\norm{X_t-X_{t-1}}^2|\mathcal{F}_{t-1}]}{\beta(t)^2} < \infty$, then $\lim\limits_{t \to \infty}  \frac{X_t}{\beta(t)}= 0$ a.s.
\label{SLNN1}
\end{lemma}

\begin{lemma}
    If $\{X_t, \mathcal{F}_t\}_{t \geq 1}$ is a non-negative submartingales, then, for any $\epsilon > 0$ we have $ \mathbb{P}(\max\limits_{1 \leq t \leq T} X(t) \geq \epsilon ) \leq \frac{E[X(T)]}{\epsilon}$.
    \label{doob's maximal inequality}
\end{lemma}
\section{Problem Statement}
Consider the following optimization problem
\begin{equation}\tag{CP1}
    \begin{split}
  \min\limits_{x \in \mathbb{X}} f(x)
    \end{split}
    \label{CP1}
\end{equation}
The constraint set $\mathbb{X}$ is a convex and compact subset of $\mathbb{R}^n$ with diameter $D$.  The function $f:  \mathbb{R}^n \to \mathbb{R}$ is a convex.   Define $ f^\ast = \min\limits_{x \in \mathbb{X}} f(x)$ and   $\mathbb{X}^\ast = \{x^\ast \in \mathbb{X} | f(x^\ast) = f^\ast\}$.
 Observe that $\mathbb{X}^\ast$ is nonempty due to compactness of the constraint set $\mathbb{X}$ and continuity of $f$. We assume  in this letter  that we have  an oracle which  generates a noisy value of the function at a given point $x \in \mathbb{X}$. That is, at each point $x\in \mathbb{X}$,  we have only the information $\hat{f}(x) = f(x) + e(x, \omega)$, where $e (x, \omega): \mathbb{R}^n \times \Omega \to \mathbb{R}$ is a random variable for each $x \in \mathbb{X}$ satisfying
\begin{equation}
\begin{split}
    & \mathbb{E}[e(x,\omega)]= b(x) \; \; \text{with} \; \; \norm{b(x)}_\ast \leq B \; \; \\ & \text{and} \; \; \mathbb{E}[\norm{e(x,\omega)}^2] \leq \mathrm{V}^2
    \end{split}
    \label{biased oracle}
\end{equation}
where, $B$ is a non-negative constant,  and $V$ can be any constant.
 \begin{remark}
  In the context of  zeroth-order  stochastic optimization problem \cite{doi:10.1137/120880811,https://doi.org/10.48550/arxiv.1809.06474,6870494}, the objective is to solve the optimization problem: $\min\limits_{x \in \mathbb{X}}  f(x) = \mathbb{E}[F(x, \omega)]$   and  the oracle only provides   $F(x, \omega)$ at any desired $x \in \mathbb{R}^n$.  In such a  situation, it is straightforward  to verify that  $\mathbb{E}[e(x,\omega)] = 0$, implying that $B = 0$. The assumption of positive $B$  makes the problem more generic than previous recent studies. In a broader sense, if $B = 0$, we call it an \emph{unbiased oracle}.
\\ However, $B$ is non-zero in many applications (see \cite{prasanth9736646} and references therein for further details), therefore the problem in this study is more general than in other recent works  due to the presence of positive $B$.
 \end{remark}

  For sake of brevity, we henceforth use $e(x)$ to denote $e(x,\omega)$.
In the next section, we discuss the zeroth-order mirror descent algorithm.
\section{Zeroth-Order Mirror Descent Algorithm}
Mirror descent algorithm is a generalization of standard subgradient descent algorithm where the Euclidean norm is replaced with a more general Bergman divergence as a proximal function. Let $R$ be the $\sigma_R$-strongly convex function and differentiable over an open set that contains the set $\mathbb{X}$. The  Bergman divergence $\mathbb{D}_R(x,y) : \mathbb{X} \times \mathbb{X} \to \mathbb{R}$ is  $\mathbb{D}_{R}(x,y) := R(x)- R(y) - \inprod{\nabla R(y)}{x-y} \; \; \forall \; \; x,y \in \mathbb{X}.$
It is clear from the definition of strong convexity that
\begin{equation}
    \mathbb{D}_R(x,y) \geq \frac{\sigma_R}{2}\norm{x-y}^2.
    \label{mirrordescent implication1}
\end{equation}
\begin{equation}
\begin{split}
    \mathbb{D}_R(z,y)- \mathbb{D}_R(z,x)-\mathbb{D}_R(x,y) &= \inprod{\nabla R(x) - \nabla R(y)}{z-x} \\
    & \forall \; \;  x,y,z \in  \mathbb{X}.
    \end{split}
    \label{mirrordescent implication2}
\end{equation}

  %
  %
 We outline  the steps of  the mirror descent algorithm.

At iteration $t$, let $x_t$ be the iterates of the ZOMD algorithm.
We  approximate the subgradient  of  function $f(x)$ at $x = x_t$ as follows.
We  generate a normal random vector $u_t \sim \mathcal{N}(0_n, I_n)$. We use the zeroth-order oracles to get the  noisy function values ($\hat{f}$) at two distinct values, that is,

$\hat{f}(x_t + \mu u_t) = f(x_t + \mu u_t) + e (x_t + \mu u_t , \omega_t^1)$ and

$\hat{f}(x_t) = f(x_t) + e (x_t,\omega_t^2)$.
Note that $\omega_t^1$ and $\omega_t^2$ are two independent realizations from the sample space $\Omega$ according to the probability law $\mathbb{P}$. Hence, we approximate the subgradient of $f$ at $x = x_t$, denoted by $\Tilde{g}(t)$ as $\Tilde{g}(t) = \frac{\hat{f}(x_t+ \mu u_t)-\hat{f}(x_t)}{\mu} u_t$.
The next iterate $x_{t+1}$ is calculated as follows:
\begin{equation}
    x_{t+1} = \argmin\limits_{x \in \mathbb{X}} \{ \inprod{\Tilde{g}(t)}{x-x_t} \} + \frac{1}{\alpha(t)} \mathbb{D}_R(x,x_t)) \}
    \label{DSMD}
\end{equation}
where, $\alpha(t)$ is the step-size of the algorithm. To show almost sure convergence, we consider weighted averaging akin to the recent work \cite{doi:10.1137/120894464} in first-order algorithm as $ z_t = \frac{\sum\limits_{j=1}^{t}\alpha(j)x_j}{\sum\limits_{k=1}^{t} \alpha(k)}$. The Bergman divergence should be chosen in such a way that \eqref{DSMD} is computationally easier to execute or a closed form solution to \eqref{DSMD} is available \cite{10.5555/1046920.1194902}.




  \begin{assumption}
The step-size $\alpha(t)$ is a decreasing sequence which satisfies  $\sum\limits_{t=1}^{\infty} \alpha(t) = \infty$ and $ \sum\limits_{t=1}^{\infty} \alpha(t)^2 < \infty.$
\label{assumption 2}
\end{assumption}

From Assumption \ref{assumption 2}, we can conclude that $\lim\limits_{t \to \infty} \alpha (t) =0$.
\begin{assumption} Let the following hold.
    \begin{enumerate}
    \item The generating random vectors $u(t) \in \mathbb{R}^n (\forall \; t \in \mathbb{N})$  are mutually independent and normally distributed and for each $t \in \mathbb{N}$ $u(t)$ is independent of $x(1),x(2),\ldots x(t)$.
        \item The random variables $e(x(t),.):\Omega \to \mathbb{R}$ and $e(x(t) + \mu u(t), .) : \Omega \to \mathbb{R}$ ( $\forall \; t \in \mathbb{N}$) are mutually independent and identically distributed in the probability space $(\Omega, \mathcal{F}, \mathbb{P})$.
        \item Define $\mathcal{H}_t = \sigma (\{x(k), u(k) \} \; 1 \leq k \leq t ) $ then  $ \mathbb{E}[e(x(t)+\mu u(t)) | \mathcal{H}_t] = b(x(t) + \mu u(t))$ and $\mathbb{E}[e(x(t)) | \mathcal{H}_t] = b(x(t))$ a.s. It is given that $\norm{b(x(t) + \mu u(t)}_\ast , \norm{b(x(t)}_\ast \leq \mathrm{B} $.  Also, $\mathbb{E}[\norm{e(x(t) + \mu u(t)}^2 | \mathcal{H}_t]$, $\mathbb{E}[\norm{e(x(t))}^2 | \mathcal{H}_t] \leq \mathrm{V}^2$ a.s.
    \end{enumerate}
    \label{independence}
\end{assumption}
Note that, most  recent literature on zeroth-order stochastic optimization computes function value  at two separate points $x(t)$ and $x(t) + \mu u_t$ under the assumption that the stochastic parameters $e(x(t))$ and $e(x(t) + \mu u_t)$ are the same. For many applications, this is rather a stringent assumption. In this letter, we avoid such an assumption, which in turn leads to significant deviation in the properties of approximated subgradient and the pertinent properties will be discussed in the ensuing section. 
        For an unbiased oracle $B =0$.
\begin{remark} Note that, most  recent literature on zeroth-order stochastic optimization computes function value  at two separate points $x_t$ and $x_t + \mu u_t$ under the assumption that the stochastic parameters $e(x_t)$ and $e(x_t + \mu u_t)$ are the same. For many applications, this is rather a stringent assumption. In this letter, we avoid such an assumption, which in turn leads to significant deviation in the properties of approximated subgradient and the pertinent properties will be discussed in the ensuing section.
\end{remark}
\section{Main Result}
In this section we discuss the properties of approximated subgradient, almost sure convergence and the finite time analysis.
 Before proceeding further, first define $\mathcal{F}_t = \sigma \{ x_l | 1\leq l \leq t\}$ $\forall \; t \in \mathbb{N}$.  Hence we get a filtration such as $\mathcal{F}_1 \subseteq \mathcal{F}_2 \subseteq \cdots \subseteq \mathcal{F}_t$. Observe that $\Tilde{g}(t-1)$ is $\mathcal{F}_t$ measurable in view of \eqref{DSMD} and also the  Bergman divergence $\mathbb{D}_R(x,x_t)$ ($\forall \; \; x \in \mathbb{X}$) is  $\mathcal{F}_t$ measurable.  Define another filtration as $\{\mathcal{G}_t\}_{t \geq 1}$ such that $\mathcal{G}_{t-1} = \mathcal{F}_t$, which will be helpful in the subsequent analysis.
\subsection{Properties of Approximated Subgradient}
The analysis in this subsection borrows some steps from \cite{NestSpok17}.
However, our analysis contains significant deviations, most notably, the result concerning the properties of approximated subgradient, which is derived using the noisy information of the  function value.



\begin{lemma}
\begin{equation*}
        \mathbb{E}[\Tilde{g}(t)|\mathcal{F}_t] = \nabla f_{\mu} (x(t)) + \mathrm{B}(t) \; \; \text{a.s.}
    \end{equation*}
    where, $\mathrm{B}(t)$ is $\mathcal{F}_t$ measurable and satisfies $\norm{\mathrm{B}(t)}_\ast \leq \frac{2 \kappa_1 \mathrm{B}}{\mu} \sqrt{n}$ a.s and we have (a.s.)
    \vspace{-0.6em}
 \begin{equation*}
    \begin{split}
        & \mathbb{E}[\norm{\Tilde{g}(t)}_\ast^2|\mathcal{F}_t]    \leq 
       \\ & \begin{cases}
            \kappa_1^2 (2 L_0^2 n + 8 \Big{(}\frac{\mathrm{V}}{\mu}\Big{)}^2 n) \; \; \text{if} \; \; f \in \mathcal{C}^{0,0}
            \\ \kappa_1^2 ( \frac{3}{4} L_1 \mu ^2 \kappa_2^4 (n+6)^3 + 3 G^2 (n+4)^2 + 12 \frac{\mathrm{V}^2}{\mu ^2}n) \; \; \text{if} \; \; f \in \mathcal{C}^{1,1} .
        \end{cases}
        \end{split}
    \end{equation*}
\label{expectation1}
\end{lemma}
\begin{proof}
It can be readily checked from the analysis in the previous section and in view of Assumption \ref{independence}  
\begin{equation}
    \mathbb{E}[\frac{f(x(t)+ \mu u(t)) -f(x(t))}{\mu} u(t) | \mathcal{F}_t] = \nabla f_\mu (x(t)).
\label{ttr}
\end{equation}
Hence, $\mathrm{B}(t) =  \mathbb{E}[\frac{e(x(t)+ \mu u(t)) -e(x(t))}{\mu} u(t) | \mathcal{F}_t]  $.  Note that because of towering property 
\begin{equation*}
\begin{split}
    & \mathbb{E}[\frac{e(x(t)+ \mu u(t)) -e(x(t))}{\mu} u(t) | \mathcal{F}_t] 
    \\ = & \mathbb{E}[\mathbb{E}[\frac{e(x(t)+ \mu u(t)) -e(x(t))}{\mu} u(t) | \mathcal{H}_t] \mathcal{F}_t]
    \end{split}
\end{equation*}
Now in view of Assumption \ref{independence}, 
\begin{flalign}
& \norm{\mathbb{E}\Big{[}\frac{e(x(t)+\mu u(t)) - e(x(t))}{\mu} u(t)| \mathcal{H}_t \Big{]}} _\ast   && \nonumber
\\ \leq  & \norm{ \frac{b(x(t) + \mu u(t)) - b(x(t))}{\mu} }_\ast\norm{u(t)}_\ast \leq \frac{2 B}{\mu} \norm{u(t)}_\ast \; \; \text{a.s.} &&
\end{flalign}
Notice that because of norm equivalence $\mathbb{E}[\norm{u(t)}_\ast | \mathcal{F}_t] \
\leq \kappa_1 \mathbb{E}[\norm{u(t)}_2|\mathcal{F}_t] \leq \kappa_1 \sqrt{n}$ a.s. Hence, $\norm{\mathrm{B}(t)}_\ast \leq \frac{2 \kappa_1 \mathrm{B}}{\mu} \sqrt{n}$ a.s.
  \\ Next we derive the second moment of the approximated subgradient. 
  Consider  the term 
  \begin{flalign}
     & \norm{\frac{f(x(t)+\mu u(t)) + e(x(t)+\mu u(t))-f(x(t)) - e(x(t))}{\mu}  u(t) }_\ast^2 && \nonumber
      \\ \leq & 2 \kappa_1^2 \norm{ \frac{f(x(t)+\mu u(t))-f(x(t))}{\mu} u(t) }_2^2  \\ + &  2 \kappa_1^2\norm{\frac{e(x(t)+\mu u(t)) - e(x(t))}{\mu}.u(t)}_2^2. && \nonumber
  \end{flalign}
 Applying the definition of $\mathcal{C}^{0,0}$, we have
 \small
 \begin{flalign}
       & \norm{\frac{f(x(t)+\mu u(t)) + e(x(t)+\mu u(t))-f(x(t)) - e(x(t))}{\mu} . u(t) }_\ast^2 && \nonumber
            \\ & \leq  2\kappa^2  L_0^2 \norm{u(t)}_2^4  + 2 \kappa_1^2 \norm{\frac{e(x(t)+\mu u(t)) - e(x(t))}{\mu}.u(t)}_2^2. &&
            \label{vari}
 \end{flalign}
 \normalsize
 \vspace{-0.6em}
    Consider the term
  \begin{flalign}
& \mathbb{E}\Big{[}\norm{\frac{e(x(t)+\mu u(t)) - e(x(t))}{\mu}u(t)}_2^2| \mathcal{H}_t \Big{]} && \nonumber
\\
\leq & \frac{2}{\mu^2} \Big{(} \mathbb{E}[(e(x(t)+\mu u(t)))^2 \norm{u(t)}_2^2 +e(x(t))^2 \norm{u(t)}_2^2|\mathcal{H}_t] \Big{)} && \nonumber
\\ \leq & \frac{4\mathrm{V}^2}{\mu^2} \norm{u(t)}^2 \; \; \text{a.s.} && \nonumber
       \end{flalign}
    Hence, by applying Towering property in \eqref{vari}
     we get,
    \begin{flalign}
        \mathbb{E}[\norm{\Tilde{g}(t)}_\ast^2|\mathcal{F}_t] \leq 2 \kappa^2 L_0^2 n + 8 \kappa_1^2 \frac{\mathrm{V}^2}{\mu^2} n \; \; \text{a.s.} && \nonumber
    \end{flalign}
For $f \in \mathcal{C}^{1,1}$, consider
\scriptsize
\begin{flalign}
    & \norm{\frac{f(x(t)+ \mu u(t)) + e(x(t)+\mu u(t))-f(x(t)) - e(x(t))}{\mu}  u(t) }_\ast^2   &&  \nonumber
      \nonumber
    \\  & \leq  3 \kappa_1^2  
 \norm{\frac{f(x(t)+\mu u(t)) -f(x(t)) - \mu \inprod{\nabla f (x(t))}{u(t)}}{\mu} u(t)}_2^2
      \\ &   + 3 \kappa_1^2 \norm{\frac{e(x(t)+\mu u(t)) - e(x(t))}{\mu}u(t)}_2^2 \nonumber
       3  \kappa_1^2 \norm{\nabla f(x(t))}_2^2\norm{u(t)}_2^4.
       \label{sml}
\end{flalign}
\normalsize
    Note that $\norm{\frac{f(x(t)+\mu u(t)) -f(x(t)) - \mu \inprod{\nabla f (x(t))}{u(t)}}{\mu} u(t)}_2^2 \leq   \frac{L_f^2 \mu ^2 \kappa_2^4}{4}   \norm{u(t)}_2^6$ because of the definition of smoothness of the function $f$.  Taking conditional expectation on $(9)$ we get the result.
\end{proof}
Using the  similar procedure we can extend the analysis  for  $f \in \mathcal{C}^{2,2}$ and so on.  It is important to note that because of consideration of more generic framework  $\mathbb{E}[\norm{\Tilde{g}(t)}_\ast^2]=  \mathcal{O}(\frac{1}{\mu^2})$ for small values of $\mu$, as opposed to  \cite{NestSpok17} because of consideration of more general framework. This result plays a significant role in the subsequent discussion of this letter.

\begin{Corollary}
    For unbiased oracle, $\mathbb{E}[\Tilde{g}_t|\mathcal{F}_t]$ $ =$ $ \nabla f_{\mu}(x_t)$ a.s.
\end{Corollary}
\subsection{Almost Sure Convergence of the ZOMD Algorithm}
Based on the discussion in Lemma \ref{expectation1},  we redefine properties of    biased subgradient as follows
\begin{equation}
      \Tilde{g}(t) = g_{\delta}(t) + \mathrm{B} (t) +\zeta (t)
        \label{bsubg}
\end{equation}
 where, $g_{\delta} (t) \in \partial_\delta f(x)$ at $x = x_t$ and  $\mathrm{B}(t)$ is $\mathcal{F}_t$  measurable and $\norm{\mathrm{B}(t)}_{\ast} \leq  B_1 $ a.s. Moreover,  $\mathbb{E}[\zeta(t)|\mathcal{F}_t] = 0$ and $ \mathbb{E}[\norm{\Tilde{g}(t)}_\ast ^2 | \mathcal{F}_t] \leq \mathrm{K}$ a.s.
Note that we can get an expression of $\delta$, $B_1$ and $K$ from Lemma \ref{expectation1}  depending on the properties of the noise and the smoothness of $f$.
\begin{theorem}
Under Assumptions \ref{assumption 2} and \ref{independence} and $\forall \; \epsilon >0$, for the iterate sequence generated by ZOMD algorithm $\{x_t\}$, there exists a subsequence $  \{x_{t_k}\}$ such that $f(x_{t_k}) - f^\ast \leq \delta + B_1D + \epsilon$ a.s.
\\ For the iterate sequence $\{z_t\}$,  $\exists \; t_0 \in \mathbb{N}$ such that  $\forall \; t \geq t_0$ we have $f(z_t)-f^\ast \leq \delta + B_1D + \epsilon$ a.s.
\label{mainths}
\end{theorem}
Before proving the Theorem \ref{mainths}, we need the following three Lemmas which we discuss here.
\begin{lemma} $\sum\limits_{t \geq 1} \frac{\alpha(t)^2}{2 \sigma_R} \norm{\Tilde{g}(t)}_\ast^2 < \infty. \; \; \text{a.s.}$
\label{lemma 6}
\end{lemma}
\begin{proof} $\lim\limits_{t \to \infty} \mathbb{E}\Big{[}\sum\limits_{k=1}^{t} \frac{\alpha(k)^2}{2 \sigma_R} \norm{\Tilde{g}(k)}_{\ast}^2 \Big{]} \leq \sum\limits_{t \geq 1} \frac{\alpha(t)^2}{2 \sigma_R} \mathrm{K} < \infty$
    By applying  Fatou's Lemma we get
    \begin{equation*}
    \mathbb{E}[\liminf\limits_{t \to \inf} \sum\limits_{k=1}^{t} \frac{\alpha(k)^2}{2\sigma_R} \norm{\Tilde{g}(k)}_\ast^2] \leq \liminf\limits_{t \to \infty} \mathbb{E}[ \sum\limits_{k=1}^{t} \frac{\alpha(k)^2}{2\sigma_R} \norm{\Tilde{g}(k)}_\ast^2]
\end{equation*}
$< \infty$.  Hence we can say $\sum\limits_{t \geq 1} \frac{\alpha(t)^2}{2\sigma_R} \norm{\Tilde{g}(t)}_\ast^2 < \infty $ a.s.
\end{proof}

\begin{lemma}
$\exists \; C >0$ such that $\mathbb{E}[\norm{\zeta(t)}_{\ast}^2| \mathcal{F}_t] < C \; \; \text{a.s.}$
\label{zeta}
\end{lemma}
\begin{proof} From the definition of $\zeta(t)$ we get that
    \begin{equation}
        \begin{split}
 \norm{\zeta(t)}_\ast^2 &  \leq 3 \kappa_1^2( \norm{\Tilde{g}(t)}_2^2 +  \norm{\mathrm{B}(t)}_2^2 +  \norm{g_{\delta}(t)}_2^2).
    \end{split}
    \label{4}
    \end{equation}
    Notice that $\exists \; K_1 >0$ such that $\norm{g_{\delta}(t)} \leq K_1$ $\forall \; t$ because of compactness of  $\mathbb{X}$.  Taking expectation on both sides of \eqref{4}, we get  (a.s.) $\mathbb{E}[\norm{\zeta(t)}_{\ast}^2|\mathcal{F}_t] \leq  3 \kappa_1^2 ( \mathrm{K} +  B_1^2 +  K_1) \delequal C. $
\end{proof}

\begin{lemma}
    \begin{equation*}
        \frac{\sum\limits_{t \geq 1}\alpha(t) \inprod{\zeta(t)}{x - x_t}}{\sum\limits_{t \geq 1} \alpha(t)} = 0 \; \; \text{a.s.} \; \; \forall \; \; x \in \mathbb{X}.
    \end{equation*}
    \label{SLLN}
\end{lemma}
\begin{proof}
    Define $X(t) = \sum\limits_{k=1}^{t} \alpha(k) \inprod{\zeta(k)}{x - x_{k}}$. In the light   of  definition of $\zeta(t)$ and since $X(t)$ is $\mathcal{F}_t$ measurable we get that $\mathbb{E}[X(t)|\mathcal{F}_t] = X(t-1)$.
    Hence $\{X(t), \mathcal{G}_t\}$ is a martingale. On the other hand, it can be seen that  (a.s.)
    \begin{equation*}
        \begin{split}
            & \sum\limits_{t \geq 1} \frac{\mathbb{E}[\norm{X(t)-X(t-1)}^2|\mathcal{F}_t]}{(\sum\limits_{k=1}^{t} \alpha(k))^2}   \leq  \\&\sum\limits_{t \geq 1} \frac{\mathbb{E}[\alpha(t)^2 \norm{\zeta(t)}_\ast^2\norm{x-x_t}^2|\mathcal{F}_t]}{(\sum\limits_{k=1}^{t} \alpha(k))^2}
             \leq \sum\limits_{t \geq 1} \frac{\alpha(t)^2 D^2C}{(\sum\limits_{k=1}^{t} \alpha(t))^2}<\infty.
        \end{split}
    \end{equation*}
    The last line is because of Lemma \ref{zeta} and the diameter of the compact set $\mathbb{X}$. Hence by applying Lemma \ref{SLNN1},  the result follows.
\end{proof}
Now we are in a position to prove the main result.

\begin{proof}
The application of first-order optimality condition to  \eqref{DSMD} yields
\begin{equation*}
\begin{split}
     & \alpha(t) \inprod{\Tilde{g}(t)}{x - x_{t+1}}
          \geq   -  \inprod{\nabla R (x_{t+1})- \nabla R (x_t)}{x-x_{t+1}}
          \end{split}
\end{equation*}
\begin{equation}
    \begin{split}
         & \geq \mathbb{D}_R(x_{t+1},x_t)+ \mathbb{D}_R(x,x_{t+1})-\mathbb{D}_R(x,x_t).
    \end{split}
    \label{firstor}
\end{equation}
The last inequality in \eqref{firstor} is due to \eqref{mirrordescent implication2}. From the LHS of \eqref{firstor}, we obtain
\begin{equation*}
    \begin{split}
        &\alpha(t) \inprod{\Tilde{g}(t)}{x - x_{t+1}} = \alpha(t) \inprod{\Tilde{g}(t)}{x-x_t+x_t-x_{t+1}}
        \\ \leq & \alpha(t) \inprod{\Tilde{g}(t)}{x-x_t} +  \frac{\alpha(t)^2}{2 \sigma_R} \norm{\Tilde{g}(t)}_\ast^2 + \frac{\sigma_R}{2} \norm{x_t-x_{t+1}}^2.
    \end{split}
\end{equation*}
The last inequality  follows by applying the Young-Fenchel inequality to the term $\alpha(t)  \inprod{\Tilde{g}(t)}{x_t - x_{t+1}}$. Hence from \eqref{firstor}, we get that
\begin{equation}
    \begin{split}
        & \mathbb{D}_R(x,x_{t+1})
        \\ \leq & \mathbb{D}_R(x,x_t) + \alpha(t) \inprod{\Tilde{g}(t)}{x-x_t} +  \frac{\alpha(t)^2}{2 \sigma_R} \norm{\Tilde{g}(t)}_\ast^2.
    \end{split}
    \label{Ber}
\end{equation}
Notice that $\mathbb{D}_R(x_{t+1},x_t) \geq \frac{\sigma_R}{2} \norm{x_{t+1}-x_t}^2$. Consider the term
\begin{equation*}
    \begin{split}
         & \alpha(t) \inprod{\Tilde{g}(t)}{x-x_t}
         =  \alpha (t) \inprod{g_{\delta} (x_t) +\mathrm{B}(t)+\zeta(t)}{x - x_t}
    \end{split}
\end{equation*}
\begin{equation}
    \begin{split}
         & \leq  \alpha(t) (f(x)-f(x_t)+ \delta + B_1D +\inprod{\zeta(t)}{x-x_t}).
    \end{split}
    \label{sps}
\end{equation}
The last inequality in \eqref{sps} is because of $\delta$-subgradient of  function $f$  and the generalized Cauchy-Schwartz inequality. Plugging  \eqref{sps}  into \eqref{Ber}
 and on applying telescopic sum  from $k = 1$ to $t$  we get,
\begin{equation}
    \begin{split}
        & \mathbb{D}_R (x^\ast,x_{t+1})
     \leq  \mathbb{D}_R(x^\ast, x_1) +  \sum\limits_{k=1}^{t} \frac{\alpha(k)^2}{2 \sigma_R} \norm{\Tilde{g}(k)}_\ast^2
    \end{split}
    \label{rumsi}
\end{equation}
\begin{equation*}
    \begin{split}
         &  +  \sum\limits_{k=1}^{t}\alpha(k) \Big{(}f^\ast-f(x_k) + \delta + B_1D  + \inprod{\zeta(k)}{x^\ast-x_k}\Big{)}.
    \end{split}
\end{equation*}
Let $\epsilon > 0$ and define the sequence of  stopping times $\{T_p\}_{p \geq 1}$ and $\{T^p\}_{p \geq 1}$ as follows:
\begin{equation*}
    \begin{split}
        T_1 & =  \inf\{f(x_t)-f^\ast \geq  \delta  + B_1D + \epsilon\}
        \\   T^1 &= \inf\{ t \geq T_1 | f(x_t) - f^\ast < \delta + B_1D + \epsilon \}
        \\ & \vdotswithin{=}
        \\  T^p &  = \inf \{ t \geq T_p | f(x_t) -f^\ast <  \delta + B_1D + \epsilon\}
        \\  T_{p+1} & = \inf\{ t \geq T^p | f(x_t)-f^\ast \geq  \delta + B_1D + \epsilon\ \}.
    \end{split}
\end{equation*}
If $\exists \; p \in \mathbb{N}$ such that infimum does not exist, we assume that $T_p = \infty$ or $T^{p} = \infty$ .

Claim-$1$ - If $T_p < \infty$, then $T^p < \infty$ a.s. $\forall \; \; p \in \mathbb{N}$.
\\ Suppose, ad absurdum,  $\exists \; p_0 \in \mathbb{N}$ such that $T_{p_0} < \infty$ but $T^{p_0} = \infty$ with probability (w.p.) $\eta$.  Let $T_{p_0} = t_{0}$, then it  implies that $\forall \; t \geq t_{0}$, $f(x_t)-f^\ast \geq \delta + B_1D + \epsilon$ w.p. $\eta$.  From \eqref{rumsi}, we deduce that $\forall \; t \geq t_0$ (w.p. $\eta$)
\begin{equation}
    \begin{split}
        & \mathbb{D}_R (x^\ast,x_{t+1}) \leq
\mathbb{D}_R(x^\ast, x_{t_0}) +      \sum\limits_{k=t_0}^{t}\alpha(k) \Big{(}-\epsilon   \\ & + \inprod{\zeta(k)}{x^\ast-x_k}\Big{)}  +  \sum\limits_{k=t_{t_0}}^{t} \frac{\alpha(k)^2}{2 \sigma_R} \norm{\Tilde{g}(k)}_\ast^2.
         \end{split}
         \label{eqn 12}
         \end{equation}
Let $t \to \infty$. Notice that in view of Lemma \ref{SLLN}, $\sum\limits_{k \geq t_{0}} \alpha(k) (- \epsilon + \inprod{\zeta(k)}{x^\ast - x_k}) = -\infty $
 and also in view of Lemma \ref{lemma 6} $\sum\limits_{k \geq t_0} \frac{\alpha(k)^2}{2 \sigma_R} \norm{\Tilde{g}(k)}_\ast^2 < \infty$ a.s.   Hence, from \eqref{eqn 12} we get
$\limsup\limits_{t \to \infty} \mathbb{D}_R(x^\ast,x_t) = - \infty$ w.p. atleast $\eta$, which implies $\eta = 0$.  Thus, $T^{p_0} < \infty$ a.s. This establishes Claim-$1$. Hence $\exists \; \{x_{t_k}\} \subseteq \{x_t\}$ such that $f(x(t_k))- f^\ast \leq \delta + B_1D + \epsilon$ a.s.

From the definition of convexity of  $f$ we get that  $\sum\limits_{k=1}^{t}\alpha(k) f (z_t) \leq \sum\limits_{j=1}^{t} \alpha(j) f(x_j)$. Hence, from \eqref{rumsi} we get that
 \begin{equation}
    \begin{split}
        & \mathbb{D}_R (x^\ast,x_{t+1})
         \leq \mathbb{D}_R(x^\ast, x_1) +    \sum\limits_{k=1}^{t} \frac{\alpha(k)^2}{2 \sigma_R} \norm{\Tilde{g}(k)}_\ast^2
          \label{rumsi1}
    \end{split}
\end{equation}
\begin{equation*}
\begin{split}
     & +\sum\limits_{k=1}^{t}\alpha(k)  \Big{(}f^\ast-f(z_t) + \delta + B_1D  + \inprod{\zeta(k)}{x^\ast-x_k}\Big{)}.
    \end{split}
\end{equation*}
In a similar fashion, define the sequence of stopping times $\{\Bar{T}_p\}_{p \geq 1}$ and $\{\Bar{T}^p\}_{p \geq 1}$ as follows:
\begin{equation*}
    \begin{split}
        & \Bar{T}_1 =  \inf\{f(z_t)-f^\ast \geq  \delta  + B_1D + \epsilon\}
        \\ &  \Bar{T}^1 = \inf\{ t \geq \Bar{T}_1 | f(z_t) - f^\ast < \delta + B_1D + \epsilon \}
        \\ & \vdotswithin{=}
        \\ & \Bar{T}^p  = \inf \{ t \geq T_p | f(z_t) -f^\ast <  \delta + B_1D + \epsilon\}
        \\ & \Bar{T}_{p+1} = \inf\{ t \geq T^p | f(z_t)-f^\ast \geq  \delta + B_1D + \epsilon\ \}.
    \end{split}
\end{equation*}
If $\Bar{T}_p < \infty$ then $\Bar{T}^p < \infty$ a.s. The reason is similar to the proof of Claim-$1$.

Claim- $2$: $\exists \; p_0 \in \mathbb{N}$ such that $\Bar{T}_{p_0} = \infty$ a.s.  If this claim is true, it proves the second part of the Theorem.

Otherwise, $\forall \;  t_1 \in \mathbb{N}$,  $\exists \; t > t_1$ such that  $(f(z_t) - f^\ast) \geq \delta + B_1D + \epsilon$ with some probability $\eta$.  Hence, from \eqref{rumsi1} we get that  \eqref{eqn 12} holds for that $t$.
         Letting $t \to \infty$ and using similar arguments  we get
 $\liminf\limits_{t \to \infty} \mathbb{D}_R(x^\ast,x_t) = -\infty$ w.p. atleast $\eta$, that means $\eta = 0$.
Hence, the Claim-$2$ holds.
\end{proof}

\begin{Corollary}[ZOMD with unbiased oracle]
     For all  $\epsilon > 0$, $\exists \; t_0 \in \mathbb{N}$ such that $\forall \; t \geq t_0$
    \begin{equation*}
        f(z_t) - f^\ast \leq
        \begin{cases}
            \mu L_0 \sqrt{n} + \epsilon \; \; \text{if} \; f \in \mathcal{C}^{0,0}
            \\ \frac{\mu ^2}{2} L_1 n + \epsilon \; \; \;  \text{if} \; f \in \mathcal{C}^{1,1} .      \end{cases}
            \; \; \text{a.s.}
    \end{equation*}
    \label{corro}
\end{Corollary}
Corollary \ref{corro} shows that by selecting a very small $\mu$, the  function value of iterate sequence converges to a small neighbourhood of the optimal value.  Notice that, $\mathbb{E}[\norm{\Tilde{g}(t)}_\ast^2] = \mathcal{O}(\frac{1}{\mu^2})$ for small value of $\mu$, hence we cannot make $\mu$ arbitrarily small. However, an analytic presentation on how small $\mu$ influences the algorithm's performance will be discussed in the ensuing section.
\begin{Corollary}[ZOMD with biased oracle]
    For all  $\epsilon > 0$ $\exists \; t_0 \in \mathbb{N}$ such that $\forall \; t \geq t_0$ the following holds
    \begin{equation*}
        f(z_t) - f^\ast \leq
        \begin{cases}
            \mu L_0 \sqrt{n} +   \frac{2 \kappa_1 B}{\mu} \sqrt{n} D+ \epsilon \; \; \text{if} \; f \in \mathcal{C}^{0,0}
            \\ \frac{\mu ^2}{2} L_1 n +  \frac{2\kappa_1 B}{\mu} \sqrt{n}D + \epsilon \; \; \;  \text{if} \; f \in \mathcal{C}^{1,1}.     \end{cases}
            \; \; \text{a.s.}
    \end{equation*}
    \label{bo}
\end{Corollary}
As Corollary \ref{bo} shows, we can not make $\mu$ very small for biased oracle. Nonetheless, an optimal $\mu^\ast$ can be calculated using Corollary \ref{bo} to show almost sure convergence to an optimal neighbourhood around the optimal value.
\subsection{Concentration Bound - Finite Time Analysis}
In the next Theorem,  we will show that a very small $\mu$ actually deteriorates the convergence rate of the ZOMD algorithm.



\begin{theorem}
Consider any  $t_0 \in \mathbb{N}$ such that $\sum\limits_{k=1}^{t_0} \alpha(k) \geq \frac{3}{\epsilon}D$. Then $\forall \; t \geq t_0$ the following holds.
    \begin{equation}
    \begin{split}
         &\mathbb{P}(f(z_t)-f^\ast \geq \delta + B_1D + \epsilon)
         \\ \leq & \frac{3\mathrm{K}}{\epsilon} \frac{\sum\limits_{k=1}^{t} \alpha(k)^2}{\sum\limits_{k=1}^{t} \alpha(k)} + \frac{9CD}{\epsilon^2} \frac{\sum\limits_{k=1}^{t} \alpha(k)^2}{(\sum\limits_{k=1}^{t} \alpha(k))^2}.
    \end{split}
    \label{theorem3}
\end{equation}
\label{concentration}
\end{theorem}

\begin{proof}
    Using the first-order optimality condition as in the proof of Theorem \ref{mainths}, we get,
    \begin{equation}
        \begin{split}
            & f(z_t)-f^\ast \leq \delta + B_1D + \frac{\mathbb{D}_R(x^\ast,x_1)}{\sum\limits_{k=1}^{t}\alpha(k)}
            \\  &
  + \frac{\sum\limits_{k=1}^{t}\alpha(k)\inprod{\zeta(k)}{x^\ast -x_k}}{\sum \limits_{k=1}^{t} \alpha(k)}
   + \frac{\sum\limits_{k=1}^{t} \alpha(k)^2\norm{\Tilde{g}(k)}_\ast^2}{2 \sigma_R \sum\limits_{k=1}^{t}\alpha(k)}.
        \end{split}
        \label{cb}
    \end{equation}
    Define $ X(t) = \sum\limits_{k=1}^{t}\alpha(k)\inprod{\zeta(k)}{x^\ast -x_k}$ and $Y(t) = \sum\limits_{k=1}^{t} \frac{\alpha(k)^2}{2 \sigma_R}\norm{\Tilde{g}(k)}_\ast^2$.
It can be seen from the definition of $\Tilde{g}(t)$ that  $ \mathbb{E}[Y(t)|\mathcal{F}_t] = Y(t-1)  + \frac{\alpha(t)^2}{2 \sigma_R} \norm{\Tilde{g}(t)}_\ast^2 \geq Y(t-1)$. Hence, $\{Y(t), \mathcal{G}_t\}$ is a non-negative sub-martingale.  It has already been shown in the proof of  Lemma \ref{SLLN} that $\{X(t), \mathcal{G}_t\}$ is  a martingale, which implies $\{\norm{X(t)}^2, \mathcal{G}_t\}$ is a sub-martingale
Choose a $t_0$ such that $\mathbb{D}_R(x^\ast,x_1) \leq \frac{\epsilon}{3} \sum\limits_{k=1}^{t} \alpha(k)$ $\forall \; t \geq t_0$ and in view of  Assumption \ref{assumption 2}, $t_0 < \infty$ .
Consider any $t > t_0$ and from \eqref{cb} if $f(z_t)-f^\ast \geq B_1D +\delta + \epsilon$ then atleast one of the following holds.
\\ $ X(t) \geq \frac{\epsilon}{3} \sum\limits_{k=1}^{t} \alpha(k)  \; \; \text{or,} \; \;  Y(t) \geq \frac{\epsilon}{3} \sum\limits_{k=1}^{t} \alpha(k)$.
That implies that $\forall \; t \geq t_0$ \begin{equation}
\begin{split}
    &\mathbb{P}(f(z_t)-f^\ast \geq \delta +
 B_1 D +   \epsilon)
    \\ \leq & \mathbb{P}(X(t) \geq \frac{\epsilon}{3} \sum\limits_{k=1}^{t} \alpha(k)) + \mathbb{P}(Y(t) \geq \frac{\epsilon}{3} \sum\limits_{k=1}^{t} \alpha(k)).
    \end{split}
    \label{ineq}
\end{equation}
Note that
\\ $\mathbb{P}(Y(t) \geq \frac{\epsilon}{3} \sum\limits_{k=1}^{t} \alpha(k)) \leq \mathbb{P}(\max\limits_{1 \leq j \leq t} Y(j) $ $\geq$ $\frac{\epsilon}{3} \sum\limits_{k=1}^{t} \alpha(k)))$. Hence, by applying Lemma \ref{doob's maximal inequality} we arrive at
\begin{equation}
    \begin{split}
        & \mathbb{P}(Y(t) \geq \frac{\epsilon}{3} \sum\limits_{k=1}^{t} \alpha(k)) \leq \frac{3}{\epsilon} \frac{E[Y(t)]}{\sum\limits_{k=1}^{t} \alpha(k)} \leq \frac{3\mathrm{K}}{\epsilon} \frac{\sum\limits_{k=1}^{t} \alpha(k)^2}{\sum\limits_{k=1}^{t} \alpha(k)}.
    \end{split}
    \label{submartingale}
\end{equation}
and similarly,
\begin{equation}
    \begin{split}
    & \mathbb{P}(X(t) \geq \frac{\epsilon}{3}\sum\limits_{k=1}^{t} \alpha(k)) \leq \mathbb{P} (\norm{X(t)}^2 \geq \frac{\epsilon^2}{9} (\sum\limits_{k=1}^{t} \alpha(k))^2) \\ &  \leq  \frac{9}{\epsilon^2} \frac{\mathbb{E}\norm{X(t)}^2}{(\sum\limits_{k=1}^{t}\alpha(k))^2} \leq \frac{9CD}{\epsilon^2} \frac{\sum\limits_{k=1}^{t} \alpha(k)^2}{(\sum\limits_{k=1}^{t} \alpha(k))^2}.
    \end{split}
    \label{martingale}
\end{equation}
Hence, by plugging \eqref{submartingale} and \eqref{martingale} into \eqref{ineq}, we get \eqref{theorem3}.
\end{proof}
\begin{remark}
     Notice that both $\mathrm{K}$ and $C$ are $\mathcal{O}(\frac{1}{\mu^2})$ from Lemma \ref{expectation}, this implies that an arbitrary small $\mu$ makes the convergence of the function value to the neighbourhood of the optimal solution slower. Hence, there is a trade-off between accuracy of the convergence to the  optimal value and convergence speed of the algorithm in the choice of $\mu$.  In the next Corollary, we capture this in detail.
\end{remark}
\begin{Corollary}
For any $\epsilon > 0$ and a confidence level $0 < p < 1$, let $p_1 = 1-p$.  Define $t_1$  such that $\forall \;  t \geq t_1$
$ \sum\limits_{k=1}^{t} \alpha(k)  \geq \frac{6 \mathrm{K}}{\epsilon p_1} \sum\limits_{k=1}^{t} \alpha(k)^2$ and $(\sum\limits_{k=1}^{t} \alpha(k))^2 \geq \frac{18 C D}{\epsilon^2 p_1} \sum\limits_{k=1}^{t} \alpha(k)^2$.
Then $\forall \; t$ $ \geq$ $ \max \{ t_0, t_1\}$ we obtain
\begin{equation*}
    \mathbb{P}( f(z_t) -f^\ast < \delta +B_1D + \epsilon ) \geq p.
\end{equation*}
Notice that  $t_1$  $< \infty$ due to Assumption \ref{assumption 2}.
\label{concentrationc}
\end{Corollary}
\section{Conclusion}
In this letter, we proved almost sure convergence of function value of ZOMD algorithm to the neighbourhood of the optimal value. Further, we derive the concentration inequality  which provides bounds on how the function value of the iterates of ZOMD algorithm deviates from the neighbourhood in any finite time. This analysis sheds some new insight to the field of zeroth-order optimization. The future path of research will attempt to demonstrate a higher convergence rate using a variance reduction technique.

\bibliographystyle{IEEEtran}
\bibliography{ref}
\end{document}